\documentclass[12pt]{article}
  \usepackage{amsfonts,amsmath,amssymb,amsthm}
  \usepackage[square,sort&compress,comma,numbers]{natbib}
  \usepackage{framed}
  \usepackage{soul}

\newcommand{\thistime}{\expandafter\calctimeA\pdfcreationdate\@nil}
\def\calctimeA#1:#2#3#4#5#6#7#8#9{\calctimeB}
\def\calctimeB#1#2#3#4#5\@nil{#1#2:#3#4}

 \usepackage{todonotes} 

   \def\R{\mathbb{R}}
   \def\N{\mathbb{N}}
   \def\Z{\mathbb{Z}}
   \def\1{{\rm I\mskip -10.5mu 1}}
   
   \def\e{{\varepsilon}}
   \def\D{{\nabla}}
   \def\lo{\mathop{\longrightarrow}}

   \def\cC{{\cal C}}

   \def\loc{\mathop{\rm loc}\nolimits}

   \newcommand{\beq}{\begin{equation}}
   \newcommand{\eeq}{\end{equation}}

\theoremstyle{definition}
\newtheorem{df}{Definition}[section]
\theoremstyle{remark}
\newtheorem{rem}[df]{Remark}
\theoremstyle{plain}
\newtheorem{prop}[df]{Proposition}
\newtheorem{lemma}[df]{Lemma}
\newtheorem{teo}[df]{Theorem}

 \newcommand{\sezione}[1]{\section{#1}\setcounter{equation}{0}}

\begin{document}


\title{On a planar Schr\"odinger-Poisson system involving a non-symmetric potential }

\author{
Riccardo Molle
\\
{\small\it Dipartimento di Matematica, Universit\`a di Roma ``Tor Vergata''}
\\
{\small\it Via della Ricerca Scientifica n. 1, 00133 Roma, Italy}
\\
{\small\tt molle@mat.uniroma2.it}
\\ \\
Andrea Sardilli
\\
{\small\tt  andreasardilli6@gmail.com}
}

\maketitle


\begin{abstract}
We prove the existence of a ground state positive solution of Schr\"odinger-Poisson systems in the plane of the form 
 $$
-\Delta u + V(x)u + \frac{\gamma}{2\pi} \left(\log|\cdot| \ast u^2 \right)u = b |u|^{p-2}u \qquad\text{in}\ \mathbb{R}^2,
$$
where $p>4$, $\gamma,b>0$  and the potential $V$ is assumed to be positive and unbounded at infinity.
On the potential we do not require any symmetry or periodicity assumption, and it is not supposed it has a limit at infinity. 
We  approach the problem by variational methods, using a variant of the mountain pass theorem and the Cerami compactness condition.
\end{abstract}


\bigskip


\sezione{Introduction} 
  
This work deals with the existence of positive solutions of Schr\"odinger-Poisson systems of the form
\begin{equation}
\tag{S}\label{P}
\left\{
\begin{array}{ll}
- \Delta u + V(x)u +  \gamma\phi(x)u = b\, |u|^{p-2}u&  \text{ in } \mathbb{R}^2\\ 
\Delta \phi = u^2&  \text{ in }  \mathbb{R}^2, 
\end{array}\right.
\end{equation}
where we assume $b>0$, the nonlinearity $p>4$, the coupling coefficient $\gamma  >0$ and the potential $V$ is a positive measurable function. 
 
For their relevance in physics, Schr\"oedinger-Poisson systems have been extensively studied in the three dimensional case, see for example the survey paper \cite{A08}, the recent paper \cite{DMT21} or also \cite{AP08,CM1,CM2,BJL13,R06} and the references therein.
On the contrary, in dimension two they are much less studied, also due to some difficulties intrinsic to the planar case   (see \cite{CW16,S08,M10,M11,DW17,CSV08,BCS21,CP21} and references therein).

Concerning the semiclassical analysis related to \eqref{P}, very little is known. 
The autonomous case has been treated first in \cite{M10}, while the case where a potential is involved has been studied only recently, in \cite{BCS21}, where the problem is fronted by perturbation methods.

We  observe that in the two dimensional Choquard problem with logarithmic kernel, analogous difficulties come out. 
We refer the reader to the recent paper \cite{CT21} for more details and references on this subject.
In \cite{CT21}, the authors considered non homogeneous nonlinearities and the positive potential is assumed to be continuous and one-periodic. 

\medskip

Since the Poisson equation in \eqref{P}  gives 
\beq
\label{ephi}
\phi=\phi_u(x) = \frac{1}{2\pi}\int_{\mathbb{R}^2} \log|x-y|\, u^2(y)\, dy,
\eeq
the system \eqref{P} can be reduced to the single nonlocal equation  
\begin{equation}
\label{E}\tag{E}
-\Delta u + V(x)u + \frac{\gamma}{2\pi} \left(\log(|\cdot|) \ast u^2 \right)u = b |u|^{p-2}u \qquad\text{in}\ \mathbb{R}^2.
\end{equation}

It is natural to expect that the solutions of equation \eqref{E} correspond to the critical points of the following action functional, where we have normalized the constants: 
$$
I(u)=\frac 12 \int_{\R^2}(|\D u|^2+V(x)u^2)dx+ \frac 14 \iint_{\R^2\times\R^2}\log(|x-y|)u^2(x)u^2(y)\, dx\, dy-\frac 1p\int_{\R^2}|u|^pdx.
$$
However, a first difficulty is that $I$ is not well defined on the usual Sobolev space $H^1(\R^2)$,  for the presence of the sign changing and unbounded logarithmic kernel.
Moreover, in the framework considered in this paper, another point to be fixed is related to the  potential, that is assumed to be unbounded and possibly not coercive at infinity.
To overcome the difficulty related to the Newtonian kernel, Stubbe in \cite{S08} introduced the following  Hilbert space, to study \eqref{P} in the autonomous case $V\equiv 0$, $b=0$,
\beq
\label{1801}
\widetilde X:=\left\{u\in H^1(\R^2)\ :\ \int_{\R^2}[\log (1+|x|)]u^2(x)\, dx<\infty\right\}. 
\eeq
The variational framework introduced in \cite{S08} was then developed in detail by Cingolani and Weth in  \cite{CW16} to front non-autonomous problems of the form \eqref{P}, with $b>0$, when the potential is assumed to be  $\Z^2$-periodic or constant. 
The Cerami condition turns out to be the more convenient compactness condition to tackle the two dimensional case in a variational way, indeed it is not evident how to show the boundedness of the Palais-Smale sequences. 
However, since the logarithmic convolution term does not have a definite sign on $\widetilde X$, also the boundedness of Cerami sequences becomes a major difficulty.
The idea, developed in \cite{S08}, to handle the different sign of the kernel is to take advantage of the relation $\log r=\log(1+r)-\log(1+\frac 1r)$, for all $r>0$, to work in $\widetilde X$ and decompose the action of the nonlocal part in its positive and negative contribution.
Once obtained the boundedness of the Cerami sequences, another difficulty is to prove that the weak limit actually is a strong limit for the Cerami sequences. 
To this aim some more assumptions on the potential has been done, for example its $\Z^2$ periodicity considered in \cite{CW16} (see also \cite{DWY21} and references therein), or its axial symmetry assumed for example in \cite{CP21,CT20}.
 
\medskip 

To the best of our knowledge, up to the semiclassical analysis performed in \cite{BCS21}, there are not works concerning problem \eqref{P} without any periodicity or symmetry assumption on the potential.
Our contribution in this work is to consider a class of problems of the form \eqref{P}, where $V$ does not enjoy of those assumptions. 
On $V$ we assume the following conditions
\beq
\tag{V}\label{HV}
\begin{array}{cl}
(a)& V\in L^1_{\loc}(\R^2)\\
(b)& \inf_{\mathbb{R}^2}V(x)>0\\ 
(c)& \left|\{x \in \R^2: V(x)\leq M \}\right| < + \infty,\ \forall\ M>0. 
\end{array}
\eeq 
If the potential verifies \eqref{HV}, the functional $I$ could be not well defined on $\widetilde X$, so that its natural domain turns out to be the weighted Hilbert space
\beq
\label{1713}
X:=\left\{u\in\widetilde X\ :\ \int_{\R^2}V(x)u^2\, dx<\infty\right\}
\eeq
endowed with the norm
$$
\|u\|_X^2:=\int_{\R^2}\left(|\D u|^2+V(x)u^2+\log(1+|x|)u^2\right)dx\qquad u\in X.
$$
We assume \eqref{HV}$(a)$ to avoid pathological potentials such that $X=\emptyset$ and to guarantee that $\cC^\infty_c(\R^2)\subseteq X$, so that the weak formulation of problem \eqref{E} is meaningful.

In Proposition \ref{PC} we employ assumption \eqref{HV}$(c)$ to get the compactness condition. 
In this regard, it is worth to observe that even if the coercivity of the logarithmic weight in the definition of $\widetilde X$ guarantees the compactness of the embedding of $\widetilde X$ in the Lebesgue spaces $L^p$, for $p\in [2,+\infty)$, the logarithmic weight in the functional $I$ does not work in the compactness of the functional, because its contribution is invariant by translation.

\medskip

We prove the following result. 
\begin{teo}
\label{T}
If $p>4$,  and $V$ verifies \eqref{HV}, 
then equation \eqref{E}  has at least a non trivial weak solution $\bar u\in X$, with $\bar u\ge 0$. 
Moreover, the solution $\bar u$ is a ground state solution:
$$
I(\bar u)=\inf\{I(u)\ :\ u\in X\setminus\{0\},\ u\mbox{ is a solution of \eqref{E}}\}.
$$
If $V\in L^q_{\loc}(\R^2)$ for some $q>1$, then $u\in\cC^{1,\alpha}_{\loc}(\R^2)$ and  $\bar u(x)>0$, for all $x\in\R^2$.

If $V\in\cC^{0,\alpha}_{\loc}(\R^2)$, then $\bar u\in \cC^2(\R^2)$ is a classical solution. 
\end{teo} 

To prove Theorem \ref{T} we use a variant of the mountain pass Theorem (see Theorem \ref{TMP}). 
Indeed it is not convenient to try a direct application of the mountain pass Theorem because it is 
not straightforward to see that $I$ is a positive quadratic form near 0 with respect to the $X$-norm. 
On the contrary, it is easy to see that $I$ is quadratic in zero with respect to the usual $H^1$-norm.
To apply Theorem \ref{TMP} a basic tool is the Cerami compactness condition, that can be verified because 
of \eqref{HV} and $p>4$.
Moreover, the assumption $p>4$ guarantees also that the fibering techniques apply and, as a consequence, the mountain pass solution turns out to be a ground state.

The paper is organized as follows: in Section 2 we introduce the main notations and some preliminary results,  
in Section 3 we prove the compactness condition and in Section 4 we prove Theorem \ref{T}.


\sezione{Variational framework and preliminary facts}

Since in the present work we focus on \eqref{E} in the case $\gamma,b,\inf V > 0$, to simplify the notations we assume
$$
\gamma = 2\pi,\qquad b=1,\qquad \inf V =1.
$$

{\small
We use the following notations:
\begin{itemize}
\item $\int=\int_{\R^2}$.
\item $B_r(y)$ denotes the open ball of radius $r>0$ and center $y\in\R^2$.
\item
$L^p$, $1\leq p<+\infty$, and $H^1$ are the usual Lebesgue and Sobolev spaces in $\R^2$, with norm
$|u|_{p}=\left(\int|u|^pdx\right)^{1/p}$ and $\|u\|=\left(\int(|\D u|^2+u^2)\,dx\right)^{1/2}$, respectively.

\item  
We denote by $C, c_1, c_2,\ldots$ various constants that can also vary from one line to another.

\end{itemize}
}

First, in Lemma \ref{Lcollect} we recall from \cite{CW16,S08} that the nonlocal term is well defined and regular in $\widetilde X$. 
To state this result, let us introduce the following bilinear forms:
$$(u, v) \mapsto B_1(u, v):= \iint \log(1+ |x-y|)u(x)v(y)\ dxdy,$$
$$(u, v)\mapsto B_2(u,v):= \iint \log\left(1+\frac{1}{|x-y|}\right)u(x) v(y)\ dxdy,$$
$$(u, v)\mapsto B_0(u, v):= B_1(u, v) - B_2(u, v) = \iint \log(|x-y|)u(x) v(y)\ dxdy,$$
where $u, v:\R^2 \to \R$ are measurable functions such that the integrals are well defined.  
Correspondingly, let 
$$
V_1(u) := B_1(u^2, u^2),\quad  V_2(u):= B_2(u^2, u^2),\quad  V_0(u):= B_0(u^2, u^2)  = V_1(u) - V_2(u).
$$
With this notations, the functional related to  problem \eqref{E} takes the form
$$I(u) = \frac12 \int (|\nabla u|^2 + V(x)u^2)\ dx + \frac14 V_0(u) - \frac 1p |u|_p^p.$$
 
\begin{lemma}
\label{Lcollect}

The functional $I$ is a well defined $C^2$ functional on the Hilbert space X. Moreover $V_2$ is of class $C^1$ on $L^\frac 83 (\R^2)$ and the following formulas hold
\begin{equation}
\label{derivateV}
V'_i(u)[v] = 4B_i(u^2, uv), \ \forall u, v \in X,\ i=0, 1, 2.
\end{equation}
In particular,
$$
V_i'(u)[u]=4V_i(u), \ \forall u  \in X,\ i=0, 1, 2.
$$
\end{lemma} 
For the proof, see  \cite[Lemma 2.2]{CW16} and \cite[Lemma 3.1.5]{Tesi_Andrea}.

 \medskip

Moreover, let us recall the Hardy-Littlewood-Sobolev inequality (see \cite[Theorem 4.3]{LL}): 
\begin{equation} \label{HLS}
\left|\iint_{\R^{N}\times \R^N} \frac {f(x) h(y)}{|x-y|^\lambda}\ dxdy\right|\leq C(\lambda, p )|f|_p|h|_{r}, \quad \forall f \in L^{p}(\R^N),\ \forall g \in L^r(\R^N),
\end{equation}
where $p,r>1$ and  $0 <\lambda< N$ are such that $\frac{1}{p} + \frac{\lambda}{2} + \frac{1}{r} = 2$.  
Then, by \eqref{HLS} and $\log(1+r)\le r$, $\forall r>0$, it is readily seen that 
 \beq
\label{stimaV_2}
V_2(u)\le \iint \frac{u^2(x)u^2(y)}{|x-y|}  \, dx\, dy\le c\, |u^2|_{4/3}  |u^2|_{4/3}=c|u|^4_{8/3}.
\eeq

\medskip

Next, we introduce the Cerami compactness condition. 

\begin{df}
Let $X$ be a Banach space,  $I\in C^1(X)$  and   $c \in \mathbb{R}$. 
A sequence $\{x_n\}_n$ is called a $(C)$ sequence for $I$ at the level $c$ if  
\begin{equation}\label{C}\tag{C}
{I(x_n) \to  c\ \text{as}\  n\to\infty \ \text{and} \ }
{\|I'(x_n)\|_{\mathcal{L}}(1 + \|x_n\|_X) \to 0,\  \text{as}\  n\to\infty.}
\end{equation}
We say that $I$ satisfies the Palais-Smale-Cerami condition (at the level $c$),  $(C)$ condition for short, if any $(C)$ sequence (at the level $c$) possesses a converging subsequence.
\end{df}

To analyze the compactness condition,  we report some known facts:  
Lemma  \ref{3.2.2} is Lemma 2.6 in \cite{CW16} and states a continuity property of the bilinear form $B_1$,
Lemma \ref{immcpt}  provides the compactness of the embedding of $X$ in the Lebesgue spaces.

\begin{lemma} 
\label{3.2.2} 
Let $\{u_n\}_n, \{v_n\}_n, \{w_n\}_n$ be bounded sequences in X such that $u_n \to u$ weakly in X. Then,  for every $z\in X$, we have $B_1(v_nw_n, z(u_n - u)) \to 0$ as $n\to \infty$.
\end{lemma}
 
\begin{lemma} 
\label{immcpt}
The space $X$ is compactly embedded in $L^q(\R^2)$ for all $q \in [2, +\infty)$.
\end{lemma}

Since bounded sets in $X$ are also bounded sets in $\widetilde X$ (see \eqref{1801}, \eqref{1713}), to verify Lemma \ref{immcpt} it is sufficient to exploit the coercivity of the map $|x|\mapsto \log(1+|x|)$  and apply, for example, Theorem XIII.65 in  \cite{RS}.

\medskip

Finally, let us state the topological tools we will use in the proof of Theorem \ref{T}. 

\medskip
\begin{lemma} [{\cite[Lemma 2.6]{LW}}]
\label{DL} 
Let $X$ be a Banach space, $I \in C^{1}(X)$, $c \in \mathbb{R}$, $\e>0$  and  suppose that there exists $\alpha >0 $ such that 
\beq
\label{1002}
(1 + \|x\|_X) \|I'(x)\|_{\mathcal{L}} \geq \alpha > 0 \qquad \forall x \in I^{c+2\epsilon}_{c-2\epsilon}.
\eeq
Then there exists a continuous map $\eta:X \to X$, such that: 
\begin{enumerate}
\item{$\eta(I^{c+\epsilon}) \subset (I^{c-\epsilon})$;}
\item{$\eta(x) = x$, for all $x \not \in I^{c+2\epsilon}_{c-2\epsilon}$},
\end{enumerate}
where $I^a=\{x\in X\ :\ I(x)\le a\}$ and $I^b_a=\{x\in X\ :\ a\le I(x)\le b\}$, for all $a<b$.
\end{lemma}

From Lemma \ref{DL} there follows the following variant of the Mountain Pass Theorem.

\begin{teo} 
\label{TMP}
Let $I \in C^{1}(X)$, X a Banach space. We suppose that $I(0) = 0$ and there exists a closed subspace $S \subset X$ that disconnects $X$ in two path-connected components $X_1$ and $X_2$. We assume that $0 \in X_1$ and there exists $A>0$, $\bar u \in X_2$ such that 
\begin{itemize}
\item{$I\big|_S\geq A$;}
\item{$I(\bar u)\leq 0$.}
\end{itemize}
Set 
\beq
\label{1011}
c = \inf_{\gamma \in \Gamma} \max_{t \in [0,1]} I(\gamma(t))
\eeq
where $\Gamma :=\{\gamma: [0,1] \to X: \gamma(0)=0,\ I(\gamma(1))\leq 0\ , \gamma(1) \in X_2\}$. 
If $I$ verifies the (C) condition at $c$, then $c$ is a critical value for $I$. 
\end{teo}

The proof of Theorem \ref{TMP} is standard so we give here only a sketch of it, 
for the sake of completeness.

If, by contradiction, $c$ is not a critical value for $I$, then \eqref{1002} has to be verified for suitable $\e,\alpha>0$, because $I$ verifies the (C) condition at the level $c$. 
We can choose $\e<A$. 
Now, let $\gamma_\e\in \Gamma$ be such that $\max_{[0,1]}I(\gamma_\e(t))\le c+\e$.
If we  consider $\tilde\gamma_\e:=\eta\circ\gamma_\e$, where $\eta$ is the deformation provided by Lemma \ref{DL}, it turns out that
$$
\tilde\gamma_\e\in\Gamma,\qquad  \max_{[0,1]}I(\tilde \gamma_\e(t))\le c-\e,
$$
contrary to the definition of $c$.


\sezione{Compactness}

This section is devoted to prove that the (C) condition holds at positive values, that is the range where we are looking for critical levels:

\begin{rem}
Let $u$ be a nontrivial critical point of $I$, then $I(u)>0$.
\end{rem}
Indeed,  by $p>4$ and assumption \eqref{HV}$(c)$, 
\begin{eqnarray*}
I(u) & = &I(u) -\frac{1}{4}I'(u)[u]  
  =  \frac{1}{4}\int \left(|\nabla u|^2dx + V(x)u^2\right)\ dx   +\left(\frac 14- \frac{1}{p}\right)|u|^p_p  \\ 
& \geq &
\frac{1}{4} \int ( |\nabla u|^2 + V(x)u^2)dx  \geq C\|u\|^2>0.
\end{eqnarray*}

\begin{prop}
\label{PC}
Assume that $V$ satisfy assumptions \eqref{HV} and let $\{u_n\}_n$ be a  $(C)$ sequence at the level $c$.
 If $c>0$ then $\{u_n\}_n$ is relatively compact.
\end{prop}

A key point to prove the proposition is the following lemma, that states the boundedness of the $(C)$ sequences at the positive levels.

\begin{lemma}
\label{LB}
Under the assumptions of Proposition \ref{PC}, the sequence $\{u_n\}_n$ is bounded in $X$.
Moreover, there exists a ball $B_2(\bar x)$ and a constant $d>0$ such that, up to a subsequence,
\beq
\label{1319}
\int_{B_2(\bar x) }u_n^2(x)\, dx \ge d \qquad \forall n\in\N.
\eeq
\end{lemma}

\proof Since $\{u_n\}_n$ is a (C) sequence, it turns out that 
\beq
\label{934}
I'(u_n)[u_n] = o(1)
\eeq
and so we have
\begin{eqnarray*}
c+o(1) &= &I(u_n) -\frac{1}{4}I'(u_n)[u_n] 
 \\ 
& = &
\frac{1}{4} \int ( |\nabla u_n|^2 + V(x)u_n^2)dx +\left(\frac 14-\frac 1p\right)\, |u_n|_p^p\\ 
& \geq &C \|u_n\|^2.
\end{eqnarray*}
Therefore,  $\{u_n\}_n$ is bounded in $H^1(\mathbb{R}^2)$,  hence in $L^q(\mathbb{R}^2)$ for all $q\in [2, +\infty)$, and  
\begin{equation}
\label{3.2}
\int_{\mathbb{R}^2} \left(|\nabla u_n|^2 + V(x)u_n^2\right)\ dx\le C\qquad\forall n\in\N.
\end{equation} 

So, we are left to prove that 
$$
|u_n|^2_{\ast}:=\int_{\R^2} \log(1 + |y|)u^2_n(x)\, dx<C,\qquad \forall n\in\N.
$$
 
Let us define
$$
\delta := \limsup_{n\to\infty} \left(\max_{y\in\mathbb{R}^2} \int_{B_1(y)}|u_n|^2dx\right).
$$ 

{\em We claim that $\delta>0$. }
Assume, by contradiction, that $\delta = 0$. 
Then by \cite[Lemma I.1]{PLL}  there follows
\beq
\label{0910}
 u_n \to 0\qquad\mbox{  in }\ L^{q}(\mathbb{R}^2),\quad\forall q\in (2,+\infty)  
\eeq
and so, by \eqref{stimaV_2}, we get
\begin{equation}
\label{0908}
V_2(u_n) \to 0 \quad \text{ as }\quad n\to+\infty.
\end{equation} 

From \eqref{934}, \eqref{0910} and \eqref{0908} we infer
$$
\int_{\mathbb{R}^2}\left(|\nabla u_n|^2 + V(x)u_n^2\right)dx + V_1(u_n) = I'(u_n)[u_n] + V_2(u_n) + |u_n|^p_p=o(1),
$$
that, by the positivity of $V_1$,  implies
\beq
\label{0945}
\int_{\mathbb{R}^2}(|\nabla u_n|^2 + V(x)u_n^2)dx=o(1)\qquad \mbox{ and }\quad
V_1(u_n)=o(1).
\eeq
Summarizing, from \eqref{0910}, \eqref{0908} and \eqref{0945} we infer
$$
I(u_n) = \frac{1}{2}\int_{\mathbb{R}^2}(|\nabla u_n|^2 + V(x)u^2_n)dx + \frac{1}{4}V_1(u_n) - \frac{1}{4} V_2(u_n) - \frac{1}{p}|u_n|^p_p =o(1),
$$
contrary to $I(u_n)\to c>0$, and the claim follows. 

\medskip

Then, $\delta > 0$ and  there exists a sequence $\{x_n\}_n \subset \mathbb{R}^2$ such that 
\beq
\label{1531}
\int_{B_1(x_n)} u^2_n(x)\,  dx > \frac{\delta}{2},
\eeq
for n large enough.  

\medskip

The task is now to prove that $\{x_n\}$ is bounded. 
Assume by contradiction that $|x_n| \to \infty$, up to a subsequence. 

Let us fix an arbitrary $M>0$ and call 
$$
A_n := \{x \in B_1(x_n)\ : \ V(x)\leq M\}, \qquad B_n := \{x \in B_1(x_n)\ : \ V(x)>M\}.
$$
By (\ref{HV})$(c)$ and $|x_n| \to \infty$ we obtain 
\beq
\label{1210}
\left|A_n\right| \leq \left|\{x \in \R^2\ :\  |x|\geq|x_n|-2, \ V(x)\leq M\}\right| \lo 0,\quad\mbox{ as }n\to\infty.
\eeq
Observe that from \eqref{1210}  we infer, for any fixed $q>2$,  
$$
\int_{A_n}u_n^2(x)\, dx\le |u|_q^2|A_n|^{1-\frac 2q}=o(1).
$$
As a consequence 
$$
\int_{B_n}u_n^2(x)\, dx  =\int_{B_1(x_n)}u_n^2(x)dx-\int_{A_n}u_n^2(x)dx\ge \frac\delta 2+o(1).
$$
Then, taking into account \eqref{stimaV_2}, we have
\beq
\label{1234}
\begin{split}
c+o(1)\geq I(u_n) &    \geq \frac{1}{2} \int_{B_n}V(x)u_n^2dx - \frac{1}{4}V_2(u_n) - C_1\frac{1}{p}\|u_n\|^p \\
& \geq  \frac\delta 4\, M - C_2,
\end{split}
\eeq
where  $C_2$ is a constant independent of $M$. 
Letting $M\to\infty$ in \eqref{1234} we get a contradiction, hence $\{x_n\}$ has to be bounded.

\medskip

Therefore there exists $\bar{x}\in \mathbb{R}^2$ such that $x_n\to \bar{x}$, up to a subsequence.    
Taking into account  \eqref{1531}, we can also assume that
\beq
\label{sem}
\int_{B_2(\bar{x})}u_n^2dx \geq \frac{\delta}{2} > 0,
\eeq
that proves \eqref{1319}.

Now, observe that for every $R>0$
$$
1 + |x - y|\geq 1 +\frac{|y|}{2}\geq \sqrt{1 + |y|} \qquad \forall y\in \mathbb{R}^2\setminus B_{2R}(0),\ \forall x\in B_R(0).
$$
Then, fixing  $\bar  R>|\bar x|+2$ and taking into account by \eqref{sem}, we get
\beq
\label{1803}
\begin{split}
V_1(u_n) =B_1(u_n^2,u_n^2) &\ge   \int_{\mathbb{R}^2\setminus B_{2\bar R}(0)}\left( \int_{B_{\bar R}(0)} \log\left(1 + |x - y|\right)u^2_n(x)\, dx\right)u^2_n(y)\, dy
 \\
&\geq  \left(\int_{B_2(\bar x)}u_n^2(x)dx\right)\cdot \left(\int_{\mathbb{R}^2\setminus B_{2\bar R}(0)} \log\left(1 + \frac{|y|}{2}\right)u^2_n(y)\, dy\right) \\
&\geq   \frac{\delta}{4}\int_{\mathbb{R}^2\setminus B_{2\bar R}(0)} \log(1 + |y|)u^2_n(y)\, dy\\
& =  \frac{\delta }{4}\left(|u_n|^2_{\ast} - \int_{B_{2\bar R}(0)} \log(1 + |y|)\, u^2_n(y)\, dy \right) \\ &\geq   \frac{\delta}{4}\left(|u_n|^2_{\ast} - \log (1+2\bar R) |u_n|^2_2  \right).
\end{split}
\eeq
Thus we conclude  that 
\begin{equation}
\label{3.10}
|u_n|^2_{\ast} \leq \frac{4}{\delta}V_1(u_n) +  C |u_n|^2_2. 
\end{equation} 
Since $\{|u_n|_2\}_n$ is bounded, if we prove that $\{V_1(u_n)\}_n$ is bounded, we are done. 

Observe that
\begin{equation}
\label{3.39}
- \frac{1}{4}V_1(u_n) + \frac{1}{4} V_2(u_n) + \left(\frac{1}{2} - \frac{1}{p}\right)|u|^p_p = I(u_n) -\frac{1}{2}I'(u_n)[u_n] = c + o(1),
\end{equation} 
from which we infer, taking into account \eqref{stimaV_2},
\begin{equation}
\label{3.4}
\frac{1}{4}V_1(u_n) = \frac{1}{4}V_2(u_n) + \left(\frac{1}{2} - \frac{1}{p}\right)|u_n|^p_p -c +o(1) 
\leq  C|u_n|^4_{\frac{8}{3}} + C|u_n|^p_p -c+o(1)\leq C,
\end{equation} 
that is the desired result.

\qed

 \medskip

{\em Proof of Proposition \ref{PC}.}\quad
By Lemma \ref{LB} $\{u_n\}_n$ is bounded in the Hilbert space $X$, so, taking also into account Lemma \ref{immcpt}, there exists $\bar u\in X$ such that, up to a subsequence, 
\beq
\label{1452}
u_n \to \bar{u}\qquad\mbox{ weakly in } X,\mbox{   in }\ L^s(\R^2), \mbox{ for every }s\in[2,+\infty), \mbox{ and a.e. in  }\R^2.
\eeq
Clearly, the sequence $\{u_n\}$ is bounded also in the weighted Hilbert space
$$
H_V:=\left\{u\in H^1(\R^2)\ :\ \int V(x)u^2(x)dx<\infty\right\},\qquad\|u\|_V^2:=\int(|\D u|^2+V(x)u^2)dx.
$$
Hence we can assume that $\{u_n\}$ weakly converges in $H_V$ to a function $\tilde u\in H_V$.
Observe that $\tilde u=\bar u$. Indeed for every fixed $\varphi\in \cC^\infty_0(\R^2)$ the map
$$
u\mapsto \int u\varphi\, dx
$$
is a continuous linear form on $X$ and on $H_V$, so that 
$$
\lim_{n\to\infty}\int u_n\varphi\, dx=\int\bar u\varphi\, dx
\quad\mbox{ and }\quad  
\lim_{n\to\infty}\int u_n\varphi\, dx=\int\tilde u\varphi\, dx.
$$
Then we can conclude that 
$$
\int (\bar u-\tilde u)\varphi\, dx =0\qquad\forall \varphi\in\cC^\infty_0(\R^2)
$$
and the assertion follows by the fundamental lemma of Calculus of Variations (see e.g. \cite[Corollary 4.24]{B}).

\medskip
 
Now, from the boundedness in $X$ of the (C) sequence $\{u_n\}_n$ we get
\beq
\label{split}
\begin{split}
o(1) = I'(u_n)[u_n - \bar{u}]   = &    \int  (|\nabla u_n|^2 + V(x)u_n^2)dx - \int \left( \nabla u_n \cdot \nabla \bar{u}+ V(x)u_n\bar{u}\right)dx
   \\
   &
+ V'_1(u_n)(u_n - \bar{u})
 - V'_2(u_n)(u_n - \bar{u})
  - \int  |u_n|^{p-2}u_n(u_n - \bar{u})\, dx.\\
  \end{split} 
  \eeq
The weak convergence in $H_V$ yields
\beq
\label{3.50}
\liminf_{n\to\infty} \int (|\nabla u_n|^2 + V(x)u_n^2)dx \geq \int  (|\nabla \bar{u}|^2 + V(x)\bar{u}^2)dx.
\end{equation}

By Lemmas \ref{Lcollect} and  \ref{3.2.2},   
  \begin{equation}
  \label{3.56}
  \begin{split} 
V'_1(u_n)[u_n - \bar{u}] & = 4B_1(u_n^2, u_n(u_n - \bar{u})) = 4B_1(u_n^2, (u_n - \bar{u})^2) + 4B_1(u_n^2, \bar{u}(u_n - \bar{u}))\\
&=4B_1(u_n^2, (u_n - \bar{u})^2)+o(1).
\end{split}
\end{equation}
By \eqref{HLS} and \eqref{1452} 
\beq
\label{1510}
\begin{split}
|V'_2(u_n)[u_n - \bar{u}]| &= |4B_2(u_n^2, u_n(u_n - \bar{u}))|  \\
& = 4\left| \iint   \log\left(1 + \frac{1}{|x - y|}\right)u^2_n(x)u_n(y)(u_n - \bar{u})(y)\,dx\,dy\right|\\
& \leq 4\iint  \frac{u^2(x)|u_n(y)|\, |(u_n - \bar{u})(y)|}{|x-y|}\, dx\,dy \\
& \leq 4  |u_n|^3_{\frac{8}{3}}|u_n - \bar{u}|_{\frac{8}{3}}=o(1).
\end{split}
\end{equation} 
 
By \eqref{1452},
\beq
\label{1942}
\left|\int  |u_n|^{p-2}u_n(u_n - \bar{u})\, dx\right|\le |u_n|_p^{p-1}|u_n-u|_p=o(1).
\eeq
 
From \eqref{split},  taking into account \eqref{3.50}, \eqref{3.56}, \eqref{1510} and \eqref{1942}, we infer
\beq
\nonumber
\begin{split}
o(1) &= \int \left(|\nabla u_n|^2 + V(x)u_n^2\right)dx - \int \left(|\nabla \bar u |^2 + V(x)\bar u^2\right)dx + 4B_1(u_n^2, (u_n - \bar{u})^2) + o(1)\\
& \geq o(1).
\end{split}
\eeq
Thus, it follows that 
\begin{equation}
\label{1813}
\int\left(|\nabla u_n|^2 + V(x)u_n^2\right)dx \to \int\left(|\nabla \bar{u}|^2 + V(x)\bar{u}^2\right)dx,
\end{equation}
and 
\begin{equation}
\label{1809}
B_1(u_n^2, (u_n - \bar{u})^2) \to 0, \ \ \text{as}\ \ n\to+\infty.
\end{equation} 
 Taking into account \eqref{1319} and arguing exactly as in \eqref{1803},  we get
\beq
\label{1810}
B_1(u_n^2, (u_n - \bar{u})^2)  \ge\frac{d}{2}(|u_n - \bar{u}|^2_*-C|u_n - \bar{u}|_2^2).
\eeq
From \eqref{1810}, \eqref{1809} and $u_n\to \bar u$ in $L^2$, there follows
$$
|u_n-\bar u|_*\to 0,
$$
that, together with \eqref{1813}, gives $\|u_n-\bar u\|_X\to 0$.
 
\qed


\sezione{Proof of Theorem \ref{T}}

To prove Theorem \ref{T}, we are going to apply Theorem \ref{TMP}. Since condition $(C)$ holds on $(0,+\infty)$ by Proposition \ref{PC}, we are left to show that also the geometric conditions hold.

Observe that, by \eqref{stimaV_2} and by the Sobolev inequality,  for every $u \in X \setminus \{0\}$  
\beq
\begin{aligned}
I(u) & \ge \frac{1}{2} \int_{\mathbb{R}^2} (|\nabla u|^2 + V(x)u^2)\ dx  - \frac{1}{4}\iint \log\left(1 + \frac{1}{|x-y|}\right)u^2(x)u^2(y)\ dxdy - \frac{1}{p}|u|^p_p\\ 
&\geq c_1 \|u\|^2 -  c_2|u|^4_{\frac{8}{3}} - c_3\|u\|^p \\
&\geq c_1\|u\|^2 - c_4\|u\|^4 - c_3\|u\|^p.
\end{aligned}
\eeq

Since $p>4$, there exist $A >0$ and $\bar{\rho}>0$ such that $I_{|_{S }}\geq A $, where  $S :=\{u \in X: \|u\| = \bar\rho \}$. 
The closed set $S$ disconnects $X$ in the two arcwise connected components 
$$
X_1:=\{u \in X: \|u\|<\bar{\rho} \}\quad\mbox{ and }\quad X_2:= \{u \in X: \|u\|>\bar{\rho} \}.
$$
Furthermore we have that $0 \in X_1$ and  there exists $\bar u\in X_2$ such that $I(\bar u)<0$, because $\lim\limits_{t\to\infty}I(tu)=-\infty$ for every $ u\in X\setminus\{0\}$.
Then, the value
\beq
\label{MPL}
c:=\inf_{\gamma \in \Gamma}\max_{t\in[0,1]} I(\gamma(t)),
\eeq
where $\Gamma := \{\gamma:[0,1] \to X: \gamma(0)=0, \gamma(1) \in X_2, I(\gamma(1))\leq0 \}$,  is a critical value for $I$, according to Theorem \ref{TMP}, with $c\ge A>0$.

To verify that $c$ is a ground state level, let us fix any nontrivial solution $\tilde u$ and consider the path $t\mapsto   t\tilde u$. 
Since $p>4$, a direct computation shows that the real function $t\mapsto I( t\tilde u)$  has a unique critical point $t_{\tilde u}$, that corresponds to a maximum. 
From $\frac{d\,}{dt}I(t\tilde u)=I'(t\tilde u)[\tilde u]$ we infer $t_{\tilde u}=1$.
Then, taking into account that $\lim\limits_{t\to\infty}I(t\tilde u)=-\infty$ and \eqref{MPL}, it is readily seen that
$$
c\le \max_{t\ge 0}I(t\tilde u)=I(\tilde u),
$$
that is our claim.

\medskip

Our next goal is to show that there exists a constant sign solution. 
Let $u$ be a critical point for $I$ at the level $c$ and consider the path $t\mapsto t|u|$. 
Since the functional is even, $c=\max_{t\ge 0}I(t|u|)=I(|u|)$ from which we infer  that $\bar u:=|u|$ is a critical point for $I$, too.
So, $\bar u$ is the weak solution we are looking for.

\medskip

Now, let us assume that $V\in L^q_{\loc}(\R^2)$ for some $q>1$ and write 
\beq
\label{1927}
-\Delta \bar u=\varphi(x)\quad\mbox{ for }\varphi(x)=-V(x)\bar u+\phi_u(x)\bar u+\bar u^{p-1},
\eeq
(see \eqref{ephi}). Since $\bar u\in L^r(\R^2)$ for every $r\in[2,\infty)$, and $\phi_u\in\cC^3(\R^2)$ (see \cite[Proposition 2.3]{CW16}), we get $\varphi\in L^s_{\loc}(\R^2)$ for every $s\in [1,q)$.
So we are in position to apply the Corollary to Theorem 2 in \cite{dB83} and obtain $\bar u\in \cC^{1,\alpha}_{\loc}(\R^2)$. Moreover, we can also write
$$
-\Delta \bar u=\psi(x)\bar u\quad\mbox{ for }\psi(x)=-V(x) +\phi_u(x) +\bar u^{p-2} 
$$
and then, since $\psi\in L^q_{\loc}(\R^2)$,  by the Harnack inequality we can conclude $\bar u>0$ (see \cite[Theorem 7.2.1]{PS07}).

If $V\in\cC^{0,\alpha}_{\loc}(\R^2)$, then the RHS in \eqref{1927} is locally H\"older continuous, so $u\in\cC^2(\R^2)$ by standard elliptic regularity (see, for example, \cite[\S 10]{LL}).

\qed


\bigskip

{\small {\bf Acknowledgement}. R.M. has been supported by the INdAM-GNAMPA group.
He acknowledges also the MIUR Excellence Department
Project awarded to the Department of Mathematics, University of Rome
Tor Vergata, CUP E83C18000100006. 
}



\end{document}